\documentclass[10pt]{article}
\usepackage{cite}
\usepackage{mathrsfs}
\usepackage{amsfonts}
\usepackage{amsmath}
\usepackage{amsfonts,amssymb}
\usepackage{dsfont}
\usepackage{curves}
\usepackage{mathrsfs}
\usepackage{pifont}
\usepackage{amssymb}
\allowdisplaybreaks

\numberwithin{equation}{section}

\date{}

\begin{document}
\title{Some sufficient conditions for path-factor uniform graphs
}
\author{\small  Sizhong Zhou$^{1}$\footnote{Corresponding
author. E-mail address: zsz\_cumt@163.com (S. Zhou)}, Zhiren Sun$^{2}$, Hongxia Liu$^{3}$\\
\small $1$. School of Science, Jiangsu University of Science and Technology,\\
\small  Zhenjiang, Jiangsu 212100, China\\
\small $2$. School of Mathematical Sciences, Nanjing Normal University,\\
\small Nanjing, Jiangsu 210023, China\\
\small  $3$. School of Mathematics and Information Sciences, Yantai University,\\
\small  Yantai, Shandong 264005, China\\
}

\maketitle
\begin{abstract}
\noindent For a set $\mathcal{H}$ of connected graphs, a spanning subgraph $H$ of $G$ is called an $\mathcal{H}$-factor of $G$ if each component
of $H$ is isomorphic to an element of $\mathcal{H}$. A graph $G$ is called an $\mathcal{H}$-factor uniform graph if for any two edges $e_1$ and
$e_2$ of $G$, $G$ has an $\mathcal{H}$-factor covering $e_1$ and excluding $e_2$. Let each component in $\mathcal{H}$ be a path with at least $d$
vertices, where $d\geq2$ is an integer. Then an $\mathcal{H}$-factor and an $\mathcal{H}$-factor uniform graph are called a $P_{\geq d}$-factor
and a $P_{\geq d}$-factor uniform graph, respectively. In this article, we verify that (\romannumeral1) a 2-edge-connected graph $G$ is a
$P_{\geq3}$-factor uniform graph if $\delta(G)>\frac{\alpha(G)+4}{2}$; (\romannumeral2) a $(k+2)$-connected graph $G$ of order $n$ with
$n\geq5k+3-\frac{3}{5\gamma-1}$ is a $P_{\geq3}$-factor uniform graph if $|N_G(A)|>\gamma(n-3k-2)+k+2$ for any independent set $A$ of $G$ with $|A|=\lfloor\gamma(2k+1)\rfloor$, where $k$ is a positive integer and $\gamma$ is a real number with $\frac{1}{3}\leq\gamma\leq1$.
\\
\begin{flushleft}

(2020) Mathematics Subject Classification: 05C70, 05C38
\end{flushleft}
\end{abstract}

\section{Introduction}

The graphs considered here are finite, undirected and simple. Let $G$ be a graph with edge set $E(G)$ and vertex set $V(G)$. We use $i(G)$,
$\omega(G)$, $\alpha(G)$ and $\delta(G)$ to denote the number of isolated vertices, the number of connected components, the independence number
and the minimum degree of $G$, respectively. Let $N_G(x)$ denote the set of neighbours of a vertex $x$ in $G$. By $d_G(x)$ we mean $|N_G(x)|$
and we call it the degree of a vertex $x$ in $G$. For any $X\subseteq V(G)$ or $X\subseteq E(G)$ the symbol $G[X]$ denotes the subgraph of
$G$ induced by $X$. We write $N_G(X)=\bigcup\limits_{x\in X}{N_{G}(x)}$ and $G-X=G[V(G)\setminus X]$ for $X\subseteq V(G)$, and denote by
$G-X$ the subgraph derived from $G$ by deleting edges of $X$ for $X\subseteq E(G)$. The edge joining vertices $x$ and $y$ is denoted by $xy$.
A vertex subset $X$ of $G$ is called an independent set if $X\cap N_G(X)=\emptyset$. Let $P_n$ and $K_n$ denote the path and the complete graph
with $n$ vertices, respectively. We denote by $K_{m,n}$ the complete bipartite graph with the bipartition $(X,Y)$, where
$|X|=m$ and $|Y|=n$.  Let $G_1$ and $G_2$ be two graphs. By $G_1\cup G_2$ we mean a graph with vertex set $V(G_1)\cup V(G_2)$ and
edge set $E(G_1)\cup E(G_2)$. By $G_1\vee G_2$ we mean a graph with vertex set $V(G_1)\cup V(G_2)$ and edge set
$E(G_1)\cup E(G_2)\cup\{e=xy: x\in V(G_1), y\in V(G_2)\}$. Recalling that $\lfloor r\rfloor$ is the greatest integer with
$\lfloor r\rfloor\leq r$, where $r$ is a real number.

A subgraph of $G$ is spanning if the subgraph includes all vertices of $G$. For a set $\mathcal{H}$ of connected graphs, a spanning subgraph
$H$ of $G$ is called an $\mathcal{H}$-factor of $G$ if each component of $H$ is isomorphic to an element of $\mathcal{H}$. A graph $G$ is
called an $\mathcal{H}$-factor covered graph if $G$ admits an $\mathcal{H}$-factor covering $e$ for any $e\in E(G)$. A graph $G$ is called an
$\mathcal{H}$-factor uniform graph if $G-e$ is an $\mathcal{H}$-factor covered graph for any $e\in E(G)$. Let each component in $\mathcal{H}$
be a path with at least $d$ vertices, where $d\geq2$ is an integer. Then an $\mathcal{H}$-factor, an $\mathcal{H}$-factor covered graph and
an $\mathcal{H}$-factor uniform graph are called a $P_{\geq d}$-factor, a $P_{\geq d}$-factor covered graph and a $P_{\geq d}$-factor uniform
graph, respectively.

Amahashi and Kano \cite{AK} derived a characterization for a graph with a $\{K_{1,l}:1\leq l\leq m\}$-factor. Kano and Saito \cite{KS} posed
a sufficient condition for the existence of $\{K_{1,l}:m\leq l\leq2m\}$-factors in graphs. Kano, Lu and Yu \cite{KLY} investigated the existence
of $\{K_{1,2},K_{1,3},K_5\}$-factors and $P_{\geq3}$-factors in graphs depending on the number of isolated vertices. Bazgan et al. \cite{BBLW}
put forward a toughness condition for a graph to have a $P_{\geq3}$-factor. Zhou, Bian and Pan \cite{ZBP}, Zhou, Sun and Liu \cite{ZSLo}, Zhou, 
Wu and Bian \cite{ZWB}, Zhou, Wu and Xu \cite{ZWX}, Zhou \cite{Zhr,Zd} obtained some results on $P_{\geq3}$-factors in graphs with given properties.
Johansson \cite{J} presented a sufficient condition for a graph to have a path-factor. Gao, Chen and Wang \cite{GCW} showed an isolated toughness
condition for the existence of $P_{\geq3}$-factors in graphs with given properties. Kano, Lee and Suzuki \cite{KLS} verified that each connected
cubic bipartite graph with at least eight vertices admits a $P_{\geq8}$-factor. Wang and Zhang \cite{WZo}, Zhou \cite{Za0}, Zhou and Liu \cite{ZL}
presented some degree conditions for the existence of graph factors. Wang and Zhang \cite{WZhr}, Yuan and Hao \cite{YH} established some relationships
between independence numbers and graph factors. Enomoto, Plummer and Saito \cite{EPS}, Zhou, Liu and Xu \cite{ZLX}, Zhou \cite{Za1,Za2}, Zhou
and Sun \cite{ZSa} derived some neighborhood conditions for the existence of graph factors. some other results on graph factors can be found in 
Wang and Zhang \cite{WZr}, Zhou, Liu and Xu \cite{ZLXb}.

A graph $H$ is factor-critical if $H-x$ has a perfect matching for each $x\in V(H)$. To characterize a graph with a $P_{\geq3}$-factor, Kaneko
\cite{K} introduced the concept of a sun. A sun is a graph formed from a factor-critical graph $H$ by adding $n$ new vertices $x_1,x_2,\cdots,x_n$
and $n$ new edges $y_1x_1,y_2x_2,\cdots,y_nx_n$, where $V(H)=\{y_1,y_2,\cdots,y_n\}$. According to Kaneko, $K_1$ and $K_2$ are also suns. A sun
with at least six vertices is called a big sun. A component of $G$ is called a sun component if it is isomorphic to a sun. Let $sun(G)$ denote
the number of sun components of $G$. Kaneko \cite{K} put forward a criterion for a graph with a $P_{\geq3}$-factor.

\medskip

\noindent{\textbf{Theorem 1.1}} (\cite{K}). A graph $G$ admits a $P_{\geq3}$-factor if and only if
$$
sun(G-X)\leq2|X|
$$
for all $X\subseteq V(G)$.

\medskip

Later, Zhou and Zhang \cite{ZZ} improved Theorem 1.1 and acquired a criterion for a $P_{\geq3}$-factor covered graph.

\medskip

\noindent{\textbf{Theorem 1.2}} (\cite{ZZ}). Let $G$ be a connected graph. Then $G$ is a $P_{\geq3}$-factor covered graph if and only if
$$
sun(G-X)\leq2|X|-\varepsilon(X)
$$
for any vertex subset $X$ of $G$, where $\varepsilon(X)$ is defined by
\[
 \varepsilon(X)=\left\{
\begin{array}{ll}
2,&if \ X \ is \ not \ an \ independent \ set;\\
1,&if \ X \ is \ a \ nonempty \ independent \ set \ and \ G-X \ has\\
&a \ non-sun \ component;\\
0,&otherwise.\\
\end{array}
\right.
\]

\medskip

Zhou and Sun \cite{ZSb} got a binding number condition for the existence of $P_{\geq3}$-factor uniform graphs. Gao and Wang \cite{GW}, Liu \cite{L} 
improved the above result on $P_{\geq3}$-factor uniform graphs. Hua \cite{H} investigated the relationship between isolated
toughness and $P_{\geq3}$-factor uniform graphs. It is natural and interesting to put forward some new sufficient conditions to guarantee that
a graph is a $P_{\geq3}$-factor uniform graph. In this article, we proceed to study $P_{\geq3}$-factor uniform graphs and pose some new graphic
parameter conditions for the existence of $P_{\geq3}$-factor uniform graphs, which are shown in the following.

\medskip

\noindent{\textbf{Theorem 1.3.}} Let $G$ be a 2-edge-connected graph. If $G$ satisfies
$$
\delta(G)>\frac{\alpha(G)+4}{2},
$$
then $G$ is a $P_{\geq3}$-factor uniform graph.

\medskip

\noindent{\textbf{Theorem 1.4.}} Let $k$ be a positive integer and $\gamma$ be a real number with $\frac{1}{3}\leq\gamma\leq1$, and let $G$ be a
$(k+2)$-connected graph of order $n$ with $n\geq5k+3-\frac{3}{5\gamma-1}$. If
$$
|N_G(A)|>\gamma(n-3k-2)+k+2
$$
for any independent set $A$ of $G$ with $|A|=\lfloor\gamma(2k+1)\rfloor$, then $G$ is a $P_{\geq3}$-factor uniform graph.

\medskip

The proofs of Theorems 1.3 and 1.4 will be given in Sections 2 and 3.

\section{The proof of Theorem 1.3}

\noindent{\it Proof of Theorem 1.3.} For any $e=xy\in E(G)$, let $G'=G-e$. To verify Theorem 1.3, we only need to prove that $G'$ is a
$P_{\geq3}$-factor covered graph. Suppose, to the contrary, that $G'$ is not a $P_{\geq3}$-factor covered graph. Then it follows from Theorem 1.2
that
\begin{align}\label{eq:2.1}
sun(G'-X)\geq2|X|-\varepsilon(X)+1
\end{align}
for some vertex subset $X$ of $G'$.

\noindent{\bf Claim 1.} $X\neq\emptyset$.

\noindent{\it Proof.} Assume that $X=\emptyset$. Then from \eqref{eq:2.1} and $\varepsilon(X)=0$ we have $sun(G')\geq1$. On the other hand, since
$G$ is 2-edge-connected, $G'$ is connected, which implies that $\omega(G')=1$. Thus, we derive that $1\leq sun(G')\leq\omega(G')=1$, that is,
$sun(G')=1$. Note that $|V(G')|=|V(G)|\geq3$ by $G$ being a 2-edge-connected graph. Hence, $G'$ is a big sun, which implies that there exist
at least three vertices $x_1,x_2,x_3$ with $d_{G'}(x_i)=1$, $i=1,2,3$. Thus, there exists at least one vertex with degree 1 in $G$, which contradicts
that $G$ is 2-edge-connected. Claim 1 is proved. \hfill $\Box$

\noindent{\bf Claim 2.} $|X|\geq2$.

\noindent{\it Proof.} Let $|X|\leq1$. Combining this with Claim 1, we get $|X|=1$.

If $G'-X$ admits a non-sun component, then $\varepsilon(X)=1$ by the definition of $\varepsilon(X)$. According to \eqref{eq:2.1} and $\varepsilon(X)=1$,
we obtain
\begin{align}\label{eq:2.2}
sun(G'-X)\geq2|X|-\varepsilon(X)+1=2|X|=2.
\end{align}

Note that $G'-X$ includes a non-sun component. Combining this with \eqref{eq:2.2}, we get
\begin{align}\label{eq:2.3}
\alpha(G')\geq sun(G'-X)+1.
\end{align}
Since $G'=G-e$, we deduce $\alpha(G)\geq\alpha(G')-2$. Then using \eqref{eq:2.2} and \eqref{eq:2.3}, we infer
\begin{align}\label{eq:2.4}
\alpha(G)\geq\alpha(G')-2\geq sun(G'-X)-1\geq2-1=1.
\end{align}

By virtue of \eqref{eq:2.2}, $G'-X$ has at least two sun components, which implies that $G-X$ admits one vertex $v$ with $d_{G-X}(v)=1$. Thus, we
derive
\begin{align}\label{eq:2.5}
\delta(G)\leq d_G(v)\leq d_{G-X}(v)+|X|=|X|+1=2.
\end{align}

It follows from \eqref{eq:2.4}, \eqref{eq:2.5} and $\delta(G)>\frac{\alpha(G)+4}{2}$ that
$$
2\geq\delta(G)>\frac{\alpha(G)+4}{2}\geq\frac{5}{2},
$$
which is a contradiction.

If $G'-X$ does not admit a non-sun component, then $\varepsilon(X)=0$ by the definition of $\varepsilon(X)$. By means of \eqref{eq:2.1}, $|X|=1$
and $\varepsilon(X)=0$, we get
\begin{align}\label{eq:2.6}
\alpha(G')\geq sun(G'-X)\geq2|X|+1=3.
\end{align}
From \eqref{eq:2.6}, we have
\begin{align}\label{eq:2.7}
\alpha(G)\geq\alpha(G')-2\geq3-2=1.
\end{align}
Note that $sun(G-X)\geq sun(G'-X)-2\geq3-2=1$ by \eqref{eq:2.6}, which implies that $G-X$ has at least one vertex $v$ with $d_{G-X}(v)\leq1$.
Thus, we infer
\begin{align}\label{eq:2.8}
\delta(G)\leq d_G(v)\leq d_{G-X}(v)+|X|\leq|X|+1=2.
\end{align}

In terms of \eqref{eq:2.7}, \eqref{eq:2.8} and $\delta(G)>\frac{\alpha(G)+4}{2}$, we derive
$$
2\geq\delta(G)>\frac{\alpha(G)+4}{2}\geq\frac{5}{2},
$$
which is a contradiction. This completes the proof of Claim 2. \hfill $\Box$

Suppose that there exist $a$ isolated vertices, $b$ $K_2$'s and $c$ big sun components $H_1,H_2,\cdots,H_c$, where $|V(H_i)|\geq6$, in $G'-X$,
and so
\begin{align}\label{eq:2.9}
sun(G'-X)=a+b+c.
\end{align}

It follows from \eqref{eq:2.1}, \eqref{eq:2.9}, $\varepsilon(X)\leq2$ and Claim 2 that
\begin{align}\label{eq:2.10}
a+b+c=sun(G'-X)\geq2|X|-\varepsilon(X)+1\geq2|X|-1\geq3.
\end{align}

\noindent{\bf Claim 3.} $\delta(G)\leq|X|+1$.

\noindent{\it Proof.} If $a\neq0$, then $d_{G'-X}(v)=0$ for any $v\in V(aK_1)$. Note that $G'=G-e$. Thus, we infer $d_{G-X}(v)\leq1$ for any
$v\in V(aK_1)$, and so
$$
\delta(G)\leq d_G(v)\leq d_{G-X}(v)+|X|\leq|X|+1.
$$

If $a=0$, then $b+c\neq0$, which implies that $G'-X$ admits at least two vertices with degree 1, and so $G-X$ has at least one vertex $v$ with
$d_{G-X}(v)=1$. Thus, we obtain
$$
\delta(G)\leq d_G(v)\leq d_{G-X}(v)+|X|=|X|+1.
$$
This completes the proof of Claim 3. \hfill $\Box$

Next, we consider two cases by the value of $a+c$.

\noindent{\bf Case 1.} $a+c=0$.

In this case, $b\geq3$ by \eqref{eq:2.10}.

\noindent{\bf Claim 4.} $\alpha(G)\geq b$.

\noindent{\it Proof.} If $x\notin V(bK_2)$ or $y\notin V(bK_2)$, then we easily see that $\alpha(G)\geq b$. If $x\in V(bK_2)$ and $y\in V(bK_2)$,
then $G-X$ has $(b-2)$ $K_2$'s and a $P_4$ component, and so we easily see that $\alpha(G)\geq(b-2)+2=b$. We finish the proof of Claim 4. \hfill $\Box$

According to \eqref{eq:2.10}, $a+c=0$, Claim 4 and $\delta(G)>\frac{\alpha(G)+4}{2}$, we deduce
\begin{align*}
\delta(G)&>\frac{\alpha(G)+4}{2}\geq\frac{b+4}{2}=\frac{a+b+c+4}{2}\\
&\geq\frac{2|X|-1+4}{2}=\frac{2|X|+3}{2}>|X|+1,
\end{align*}
which contradicts Claim 3.

\noindent{\bf Case 2.} $a+c\neq0$.

\noindent{\bf Subcase 2.1.} $a\neq0$.

If $x\notin V(aK_1)$ and $y\notin V(aK_1)$, then $d_{G-X}(v)=0$ for any $v\in V(aK_1)$. Thus, we derive
\begin{align}\label{eq:2.11}
\delta(G)\leq d_G(v)\leq d_{G-X}(v)+|X|=|X|.
\end{align}
It follows from \eqref{eq:2.10}, \eqref{eq:2.11}, $\delta(G)>\frac{\alpha(G)+4}{2}$ and $\alpha(G)\geq sun(G-X)\geq sun(G'-X)-2$ that
\begin{align*}
|X|&\geq\delta(G)>\frac{\alpha(G)+4}{2}\geq\frac{sun(G'-X)-2+4}{2}=\frac{sun(G'-X)+2}{2}\\
&\geq\frac{2|X|-1+2}{2}=|X|+\frac{1}{2},
\end{align*}
which is a contradiction. In what follows, we discuss the case with $x\in V(aK_1)$ or $y\in V(aK_1)$. Without loss of generality, let $x\in V(aK_1)$.
We write $Y=V(H_1)\cup\cdots\cup V(H_c)$.

\noindent{\bf Subcase 2.1.1.} $y\in V(bK_2)\cup Y$.

In this subcase, we deduce $\alpha(G)\geq a+b+c$. Combining this with \eqref{eq:2.10} and $\delta(G)>\frac{\alpha(G)+4}{2}$, we infer
$$
\delta(G)>\frac{\alpha(G)+4}{2}\geq\frac{a+b+c+4}{2}\geq\frac{2|X|-1+4}{2}=\frac{2|X|+3}{2}>|X|+1,
$$
which contradicts Claim 3.

\noindent{\bf Subcase 2.1.2.} $y\in V(G)\setminus(V(bK_2)\cup Y)$.

In this subcase, we have $sun(G-X)\geq sun(G'-X)-1$. Combining this with \eqref{eq:2.10}, $\alpha(G)\geq sun(G-X)$ and $\delta(G)>\frac{\alpha(G)+4}{2}$,
we derive
\begin{align*}
\delta(G)&>\frac{\alpha(G)+4}{2}\geq\frac{sun(G-X)+4}{2}\geq\frac{sun(G'-X)-1+4}{2}\\
&=\frac{sun(G'-X)+3}{2}\geq\frac{2|X|-1+3}{2}=|X|+1,
\end{align*}
which contradicts Claim 3.

\noindent{\bf Subcase 2.2.} $c\neq0$.

Obviously, $\alpha(G')\geq a+b+\sum\limits_{i=1}^{c}{\frac{|V(H_i)|}{2}}\geq a+b+3c$ by $|V(H_i)|\geq6$. Combining this with \eqref{eq:2.10}, $c\neq0$
and $\alpha(G)\geq\alpha(G')-2$, we admit
\begin{align}\label{eq:2.12}
\alpha(G)\geq\alpha(G')-2\geq a+b+3c-2\geq a+b+c\geq2|X|-1.
\end{align}

By virtue of \eqref{eq:2.12}, Claim 3 and $\delta(G)>\frac{\alpha(G)+4}{2}$, we deduce
$$
|X|+1\geq\delta(G)>\frac{\alpha(G)+4}{2}\geq\frac{2|X|-1+4}{2}=|X|+\frac{3}{2},
$$
which is a contradiction. This completes the proof of Theorem 1.3. \hfill $\Box$

\section{The proof of Theorem 1.4}

\noindent{\it Proof of Theorem 1.4.} For any $e\in E(G)$, we write $G'=G-e$. To prove Theorem 1.4, we only need to justify that $G'$ is a
$P_{\geq3}$-factor covered graph. Suppose, to the contrary, that $G'$ is not a $P_{\geq3}$-factor covered graph. Then by Theorem 1.2, we admit
\begin{align}\label{eq:3.1}
sun(G'-X)\geq2|X|-\varepsilon(X)+1
\end{align}
for some vertex subset $X$ of $G'$. We write $a=i(G-X)$ and $b=\lfloor\gamma(2k+1)\rfloor$.

\noindent{\bf Claim 1.} $b\geq a+1$.

\noindent{\it Proof.} Let $b\leq a$. We may choose $b$ isolated vertices $x_1,x_2,\cdots,x_b$ in $G-X$. Write $A=\{x_1,x_2,\cdots,x_b\}$. Then
$A$ is an independent set of $G$. Thus, we infer
\begin{align}\label{eq:3.2}
\gamma(n-3k-2)+k+2<|N_G(A)|\leq|X|.
\end{align}

It follows from \eqref{eq:3.1}, \eqref{eq:3.2} and $\varepsilon(X)\leq2$, $\frac{1}{3}\leq\gamma\leq1$ and $n\geq5k+3-\frac{3}{5\gamma-1}$ that
\begin{align*}
0&\geq|X|+sun(G'-X)-n\geq|X|+2|X|-\varepsilon(X)+1-n\\
&\geq3|X|-n-1>3(\gamma(n-3k-2)+k+2)-n-1\\
&=(3\gamma-1)n-3\gamma(3k+2)+3k+5\\
&\geq(3\gamma-1)\left(5k+3-\frac{3}{5\gamma-1}\right)-3\gamma(3k+2)+3k+5\\
&=(3\gamma-1)(2k+1)-\frac{3(3\gamma-1)}{5\gamma-1}+3\\
&\geq3-\frac{3(3\gamma-1)}{5\gamma-1}>3-3=0,
\end{align*}
which is a contradiction. We finish the proof of Claim 1. \hfill $\Box$

In what follows, we consider four cases by the value of $|X|$ and derive a contradiction in each case.

\noindent{\bf Case 1.} $|X|=0$.

Note that $G'=G-e$ and $G$ is $(k+2)$-connected. Hence, $G'$ is $(k+1)$-connected and $\omega(G')=1$. Combining this with \eqref{eq:3.1} and
$\varepsilon(X)=0$, we obtain $1=\omega(G')\geq sun(G')\geq1$. Thus, we have $sun(G')=\omega(G')=1$. Then using
$n\geq5k+3-\frac{3}{5\gamma-1}\geq8-\frac{3}{5\times\frac{1}{3}-1}=\frac{7}{2}>3$, we see that $G'$ is a big sun, and so $G'$ has at least
three vertices with degree 1, which contradicts that $G'$ is a $(k+1)$-connected graph.

\noindent{\bf Case 2.} $1\leq|X|\leq k$.

Note that $1\leq|X|\leq k$ and $G'$ is $(k+1)$-connected. We derive $\omega(G'-X)=1$. According to \eqref{eq:3.1} and $\varepsilon(X)\leq|X|$,
we get
$$
1=\omega(G'-X)\geq sun(G'-X)\geq2|X|-\varepsilon(X)+1\geq|X|+1\geq2,
$$
which is a contradiction.

\noindent{\bf Case 3.} $|X|=k+1$.

Since $G$ is $(k+2)$-connected, $G-X$ is connected, and so $\omega(G-X)=1$. Note that $G'=G-e$. Thus, we deduce
\begin{align}\label{eq:3.3}
\omega(G'-X)\leq\omega(G-X)+1=2.
\end{align}
By virtue of \eqref{eq:3.1}, \eqref{eq:3.3}, $k\geq1$ and $\varepsilon(X)\leq2$, we infer
\begin{align*}
2&\geq\omega(G'-X)\geq sun(G'-X)\geq2|X|-\varepsilon(X)+1\geq2|X|-1\\
&=2(k+1)-1=2k+1\geq3,
\end{align*}
which is a contradiction.

\noindent{\bf Case 4.} $|X|\geq k+2$.

In light of \eqref{eq:3.1}, $\varepsilon(X)\leq2$ and $\frac{1}{3}\leq\gamma\leq1$, we derive
\begin{align*}
sun(G-X)&\geq sun(G'-X)-2\geq2|X|-\varepsilon(X)+1-2\geq2|X|-3\\
&\geq2(k+2)-3=2k+1\geq\gamma(2k+1)\geq\lfloor\gamma(2k+1)\rfloor=b,
\end{align*}
which implies that $G-X$ admits an independent set of order at least $b$. Then using Claim 1, we may choose $a$ isolated vertices $x_1,x_2,\cdots,x_a$
and $(b-a)$ nonadjacent vertices $x_{a+1},\cdots,x_b$ with $d_{G-X}(x_i)=1$ for $a+1\leq i\leq b$, in $G-X$. Set $A=\{x_1,x_2,\cdots,x_a,x_{a+1},\cdots,x_b\}$.
Then $A$ is an independent set of $G$. Thus, we deduce
$$
\gamma(n-3k-2)+k+2<|N_G(A)|\leq|X|+b-a,
$$
that is,
\begin{align}\label{eq:3.4}
|X|>\gamma(n-3k-2)+k+2-b+a.
\end{align}
It follows from \eqref{eq:3.1}, \eqref{eq:3.4}, $\varepsilon(X)\leq2$ and $n\geq5k+3-\frac{3}{5\gamma-1}$ that
\begin{align*}
0&\geq|X|+2sun(G'-X)-i(G'-X)-n\\
&\geq|X|+2(2|X|-\varepsilon(X)+1)-(i(G-X)+2)-n\\
&\geq|X|+2(2|X|-1)-(a+2)-n\\
&=5|X|-a-4-n\\
&>5(\gamma(n-3k-2)+k+2-b+a)-a-4-n\\
&=(5\gamma-1)n-5\gamma(3k+2)+5k+10-5b+4a-4\\
&\geq(5\gamma-1)\left(5k+3-\frac{3}{5\gamma-1}\right)-5\gamma(3k+2)+5k+6-5b\\
&=5\gamma(2k+1)-5b\\
&=5\gamma(2k+1)-5\lfloor\gamma(2k+1)\rfloor\\
&\geq0,
\end{align*}
which is a contradiction. This completes the proof of Theorem 1.4. \hfill $\Box$

\section{Remarks}

\noindent{\bf Remark 1.} Next, we show that the condition $\delta(G)>\frac{\alpha(G)+4}{2}$ in Theorem 1.3 cannot be replaced by
$\delta(G)\geq\frac{\alpha(G)+4}{2}$. We construct a graph $G=K_{3+t}\vee(4+2t)K_2$, where $t$ is an nonnegative integer. Then $G$ is $(3+t)$-connected,
$\delta(G)=4+t$ and $\alpha(G)=4+2t$. Thus, we admit $\delta(G)=\frac{\alpha(G)+4}{2}$. For any $e\in E((4+2t)K_2)$, let
$G'=G-e=K_{3+t}\vee((3+2t)K_2\cup(2K_1))$. Select $X=V(K_{3+t})\subseteq V(G')$. Then $|X|=3+t$ and $\varepsilon(X)=2$. Thus, we derive
$$
sun(G'-X)=5+2t>4+2t=2(3+t)-2=2|X|-\varepsilon(X).
$$
By means of Theorem 1.2, $G'$ is not a $P_{\geq3}$-factor covered graph, and so $G$ is not a $P_{\geq3}$-factor uniform graph.

\medskip

\noindent{\bf Remark 2.} The conditions with a $(k+2)$-connected graph and $|N_G(A)|>\gamma(n-3k-2)+k+2$ in Theorem 1.4 cannot be replaced by
a $(k+1)$-connected graph and $|N_G(A)|\geq\gamma(n-3k-2)+k+1$. Let $\gamma$ be a rational number such that $\frac{1}{3}\leq\gamma\leq1$. Then
we can write $\gamma=\frac{b}{2k+1}$ for nonnegative integers $b$ and $k$. Let $G=K_{k+1}\vee((2k+1)K_2)$, where $k$ is a positive integer.
Then $G$ is $(k+1)$-connected and $n=|V(G)|=5k+3>5k+3-\frac{3}{5\gamma-1}$. If $A$ is an independent set of order $b=\gamma(2k+1)$, then
$$
\gamma(n-3k-2)+k+2>|N_G(A)|=\gamma(2k+1)+k+1=\gamma(n-3k-2)+k+1.
$$
For any $e\in E((2k+1)K_2)$, let $G'=G-e=K_{k+1}\vee((2k)K_2\cup(2K_1))$. Select $X=V(K_{k+1})\subseteq V(G')$. Then $|X|=k+1$ and $\varepsilon(X)=2$.
Thus, we infer
$$
sun(G'-X)=2k+2>2k=2(k+1)-2=2|X|-\varepsilon(X).
$$
According to Theorem 1.2, $G'$ is not a $P_{\geq3}$-factor covered graph, and so $G$ is not a $P_{\geq3}$-factor uniform graph.

\section{Conclusion}

The concept of path-factor uniform graph was first presented by Zhou and Sun \cite{ZSb}, and they showed a binding number condition for the existence
of $P_{\geq3}$-factor uniform graphs. Gao and Wang \cite{GW}, Liu \cite{L} improved Zhou and Sun's above result. Hua \cite{H} gave toughness and isolated toughness
conditions for graphs to be $P_{\geq3}$-factor uniform graphs. In our article, we study the relationships between some graphic parameters (for instance,
minimum degree, independence number and neighborhood, and so on) and the existence of $P_{\geq3}$-factor uniform graphs. The theorems derived in this
article belong to existence theorems, that is, under what kind of conditions the path-factor uniform graph exists. However, in a specific computer
network, it needs to use a certain algorithm to determine the values of some graphic parameters of the fix network graph and show the eligible
path-factor uniform graph from the algorithm point of view. The problems of such algorithms are worthy of consideration in the future research.

Until now, the results on the existence of path-factor uniform graphs are very few. There are many problems on graphs which are path-factor uniform
graphs can be considered. For example, we can consider the structures and properties of path-factor uniform graphs. In what follows, we put forward
the open problems as the end of our article.

\medskip

\noindent{\bf Problem 1.} Find the necessary and sufficient conditions for a graph to be a path-factor uniform graph.

\medskip

\noindent{\bf Problem 2.} Find the relationships between the other graphic parameters and path-factor uniform graphs.

\medskip

\noindent{\bf Problem 3.} What are the structures and properties in path-factor uniform graphs?

\medskip



\begin{thebibliography}{9999}

\bibitem {AK} A. Amahashi, M. Kano, On factors with given components, Discrete Mathematics 42(1982)1--6.

\bibitem {KS} M. Kano, A. Saito, Star-factors with large component, Discrete Mathematics 312(2012)2005--2008.

\bibitem {KLY} M. Kano, H. Lu, Q. Yu, Component factors with large components in graphs, Applied Mathematics Letters 23(2010)385--389.

\bibitem {BBLW} C. Bazgan, A. Benhamdine, H. Li, M. Wo\'zniak, Partitioning vertices of 1-tough graph into paths, Theoretical Computer Science 263(2001)255--261.

\bibitem {ZBP} S. Zhou, Q. Bian, Q. Pan, Path factors in subgraphs, Discrete Applied Mathematics, DOI: 10.1016/j.dam.2021.04.012

\bibitem {ZSLo} S. Zhou, Z. Sun, H. Liu, On $P_{\geq3}$-factor deleted graphs, Acta Mathematicae Applicatae Sinica-English Series 38(1)(2022)178--186.

\bibitem {ZWB} S. Zhou, J. Wu, Q. Bian, On path-factor critical deleted (or covered) graphs, Aequationes Mathematicae, DOI: 10.1007/s00010-021-00852-4

\bibitem {ZWX} S. Zhou, J. Wu, Y. Xu, Toughness, isolated toughness and path factors in graphs, Bulletin of the Australian Mathematical Society,
DOI: 10.1017/S0004972721000952

\bibitem {Zhr} S. Zhou, Remarks on path factors in graphs, RAIRO-Operations Research 54(6)(2020)1827--1834.

\bibitem{Zd} S. Zhou, Degree conditions and path factors with inclusion or exclusion properties, Bulletin Mathematique de la Societe des Sciences
Mathematiques de Roumanie, Accept.

\bibitem {J} R. Johansson, An El-Zah\'ar type condition ensuring path-factors, Journal of Graph Theory 28(1998)39--42.

\bibitem {GCW} W. Gao, Y. Chen, Y. Wang, Network vulnerability parameter and results on two surfaces, International
Journal of Intelligent Systems 36(2021)4392--4414.

\bibitem {KLS} M. Kano, C. Lee, K. Suzuki, Path and cycle factors of cubic bipartite graphs, Discussiones Mathematicae Graph Theory 28(3)(2008)551--556.

\bibitem {WZo} S. Wang, W. Zhang, On $k$-orthogonal factorizations in networks, RAIRO-Operations Research 55(2)(2021)969--977.

\bibitem {Za0} S. Zhou, A Fan-type result for the existence of restricted fractional $(g,f)$-factors, Proceedings of the Romanian Academy, Series A:
Mathematics, Physics, Technical Sciences, Information Science 22(1)(2021)3--10.

\bibitem {ZL} S. Zhou, H. Liu, Discussions on orthogonal factorizations in digraphs, Acta Mathematicae Applicatae Sinica-English Series 38(2)(2022)417--425.

\bibitem {WZhr} S. Wang, W. Zhang, Remarks on fractional ID-$[a,b]$-factor-critical covered network graphs, Proceedings of the Romanian Academy,
Series A: Mathematics, Physics, Technical Sciences, Information Science 22(3)(2021)209--216.

\bibitem {YH} Y. Yuan, R. Hao, Independence number, connectivity and all fractional $(a,b,k)$-critical graphs, Discussiones Mathematicae Graph
Theory 39(2019)183--190.

\bibitem {EPS} H. Enomoto, M. Plummer, A. Saito, Neighborhood unions and factor critical graphs, Discrete Mathematics 205(1999)217--220.

\bibitem {ZLX} S. Zhou, H. Liu, Y. Xu, A note on fractional ID-$[a,b]$-factor-critical covered graphs, Discrete Applied Mathematics,
DOI: 10.1016/j.dam.2021.03.004

\bibitem{Za1} S. Zhou, A neighborhood union condition for fractional $(a,b,k)$-critical covered graphs, Discrete Applied Mathematics,
DOI: 10.1016/j.dam.2021.05.022

\bibitem {Za2} S. Zhou, A result on fractional $(a,b,k)$-critical covered graphs, Acta Mathematicae Applicatae Sinica-English Series 37(4)(2021)657--664.

\bibitem {ZSa} S. Zhou, Z. Sun, A neighborhood condition for graphs to have restricted fractional $(g,f)$-factors, Contributions to Discrete
Mathematics 16(1)(2021)138--149.

\bibitem {WZr} S. Wang, W. Zhang, Research on fractional critical covered graphs, Problems of Information Transmission 56(2020)270--277.

\bibitem {ZLXb} S. Zhou, H. Liu, Y. Xu, Binding numbers for fractional $(a,b,k)$-critical covered graphs, Proceedings of the Romanian Academy,
Series A: Mathematics, Physics, Technical Sciences, Information Science 21(2)(2020)115--121.

\bibitem {K} A. Kaneko, A necessary and sufficient condition for the existence of a path factor every component of which is a path of length at
least two, Journal of Combinatorial Theory, Series B 88(2003)195--218.

\bibitem {ZZ} H. Zhang, S. Zhou, Characterizations for $P_{\geq2}$-factor and $P_{\geq3}$-factor covered graphs, Discrete Mathematics
309(2009)2067--2076.

\bibitem {ZSb} S. Zhou, Z. Sun, Binding number conditions for $P_{\geq2}$-factor and $P_{\geq3}$-factor uniform graphs, Discrete Mathematics 343(3)(2020)111715.

\bibitem {GW} W. Gao, W. Wang, Tight binding number bound for $P_{\geq3}$-factor uniform graphs, Information Processing Letters 172(2021)106162.

\bibitem {L} H. Liu, Binding number for path-factor uniform graphs, Proceedings of the Romanian Academy, Series A: Mathematics, Physics, Technical
Sciences, Information Science 23(1)(2022)25--32.

\bibitem {H} H. Hua, Toughness and isolated toughness conditions for $P_{\geq3}$-factor uniform graphs, Journal of Applied Mathematics and Computing
66(2021)809--821.


\end{thebibliography}
\end{document}